\input amstex
\documentstyle{amsppt}
\magnification=1200
\NoRunningHeads
\TagsOnRight

\define\R{{\Bbb R}}

\define\Z{{\Bbb Z}}

\define\E{{\Bbb E}}

\define\wtQ{\widetilde Q}
\define\wtrho{\widetilde\rho}
\define\whQ{\widehat Q}
\define\whk{\widehat k}
\define\Om{\Omega}
\define\om{\omega}
\define\al{\alpha}
\define\be{\beta}

\define\Ga{\Gamma}
\define\de{\delta}

\define\ka{\kappa}
\define\la{\lambda}

\define\si{\sigma}
\define\ep{\varepsilon}

\define\Pz{P_{zz'}}
\define\PP{{\Cal P}}
\define\PPz{\PP_{zz'}}
\define\PDt{\Cal{PD}(t)}

\define\rz{\rho^{(zz')}}

\define\tht{\thetag}

\define\sh{\operatorname{sh}}
\define\ch{\operatorname{ch}}
\define\tg{\operatorname{tg}}
\define\th{\operatorname{th}}
\define\Det{\operatorname{Det}}
\define\Prob{\operatorname{Prob}}
\define\const{\operatorname{const}}
\define\Nt{N_\tau}

\topmatter
\title Point processes and\\ the infinite symmetric group \\
Part III: Fermion point processes
\endtitle

\author Alexei Borodin and Grigori Olshanski
\endauthor

\thanks Supported by the Russian Foundation for Basic Research under
grant 98-01-00303 (G.~O)
and by the Russian Program for Support of Scientific Schools under
grant 96-15-96060 (A.~B. and G.~O.)
\endthanks
\abstract
In Part I (G.~Olshanski) and Part II (A.~Borodin) we developed an
approach to certain probability distributions on the Thoma simplex.
The latter has infinite dimension and is a kind of dual object for the
infinite symmetric group.  Our approach is based on studying the
correlation functions of certain related point stochastic processes.

In the present paper we consider the so--called tail point processes
which describe the limit behavior of the Thoma parameters
(coordinates on the Thoma simplex) with large numbers. The tail
processes turn out to be stationary processes on the real line. Their
correlation functions have determinantal form with a kernel which
generalizes the well--known sine kernel arising in random matrix theory.
Our second result is a law of large numbers for the Thoma parameters.
We also produce Sturm--Liouville operators commuting with the Whittaker
kernel introduced in Part II and with the generalized sine kernel.
\endabstract

\toc
\widestnumber\head{1}
\head {} Introduction \endhead
\head 1. The fermion point processes \endhead
\head 2. Lifting and the Whittaker kernel \endhead
\head 3. The tail process \endhead
\head 4. The $\sin/\sh$ and $\sh/\sh$ kernels \endhead
\head 5. Rate of decay of the Thoma parameters \endhead
\head 6. Associated Sturm--Liouville operators \endhead
\head 7. Comparison with Poisson--Dirichlet \endhead
\head {} References \endhead
\endtoc

\endtopmatter

\document

The present note is a continuation of \cite{O} (Part I) and \cite{B1}
(Part II). Our aim here is to draw some conclusions from the
computations of Part II.

In Part I we started the study of a family $\{\Pz\}$ of probability
Borel measures living on an infinite--dimensional simplex $\Om$, the
Thoma simplex. Recall that this family consists of two parts: the
principal series and the complementary series, each of which is
indexed by two real parameters. For the principal
series, $z$ is an arbitrary complex number distinct from
$0,\pm1,\pm2,\dots$, and $z'=\bar z$. For the complementary series,
$z$ and $z'$ are real numbers which are both contained in a unit
interval with integer ends. When $z$ is a noninteger real number and
$z'=z$, the corresponding measure belongs to the intersection of the
both series. 

Note that $z-z'$ ranges over the imaginary axis plus the open
interval $(-1,1)$ -- a picture that immediately evokes the
principal and complementary series for $SL(2,\R)$.  

The measures $\Pz$ originated from the work \cite{KOV}: they govern
the decomposition of certain reducible representations which seem to be
right analogs of the regular representation for the infinite
symmetric group. It is known that, excepting the symmetry relation
$\Pz=P_{z'z}$, the measures $\Pz$ are pairwise disjoint. 

Recall that the points of the Thoma simplex $\Om$ are the double
sequences $\al=(\al\ge\al_2\ge\dots\ge0)$, 
$\be=(\be_1\ge\be_2\ge\dots\ge0)$ such that $\sum(\al_j+\be_j)\le1$.
The basic idea of Part I was to interpret the measures $\Pz$ as
point processes (denoted as $\PPz$) on the punctured interval
$I=[-1,1]\setminus\{0\}$, the random configuration being of the form
$(-\be_1,-\be_2,\dots,\al_2,\al_1)$, where only nonzero $\al_j$'s and
$\be_j$'s are considered and the points are accumulated near 0. A
method of calculating the correlation functions $\rz_n$,
$n=1,2,\dots$, of the processes $\PPz$ was developed. In Part I, we
calculated the first correlation function $\rz_1$, and in Part II ---
the higher correlation functions $\rz_n$. More advanced results were
obtained for the process $\PPz^+$, the restriction of $\PPz$ to
$(0,1]\subset I$. The process $\PPz^+$ reflects the behavior of the
Thoma parameters $\al_j$; the study of the $\be_j$'s is reduced to
that of the $\al_j$'s simply by change of the sign for $z$ and $z'$.

The initial definition of the processes $\PPz$ is rather indirect 
and we know no explicit probabilistic mechanism generating them.
{}From the beginning it was unclear what known processes they
resemble. Now, the knowledge of the correlation functions makes it
possible to conclude that $\PPz$ (or at least certain derived
processes) are similar to the point processes arising in the scaling
limit of certain random matrix ensembles. The basic common feature is
that the correlation functions are given by determinantal formulas
involving a kernel; such processes are called fermion point
processes after \cite{Ma1, Ma2, DV}. Note that determinantal correlation
functions also appear in certain models of mathematical physics \cite{KBI}. 

The results of Part II lead to interesting kernels: the Whittaker
kernel (see sections 1--2), the $\sin/\sh$ kernel, the $\sh/\sh$
kernels and their degenerations (see sections 3--4). The Whittaker
kernel seems to be a new example; the $\sin/\sh$ kernel already appeared 
in works of mathematical physicists, see \cite{BCM, MCIN}. We think
that the connection of our problem with the random matrix theory is
interesting and promising.  

The note is organized as follows. In section 1, we briefly review
some general facts about the fermion processes. In section 2, we
discuss the ``lifting'' of the process $\PPz^+$, which leads to the
Whittaker kernel. As an application, we calculate the mean value for
$\sum\al_i$ and $\sum\be_i$. In section 3, we introduce the ``tail
process'' for $\PPz^+$; it turns out to be a stationary fermion
process on $\R$ depending on $z,z'$. In this way, we get a
two--parametric family of kernels generalizing the sine
kernel; they are discussed in section 4. In section 5, we prove that
$$
\lim_{j\to\infty}\al_j^{1/j}=\lim_{j\to\infty}\be_j^{1/j}=e^{-1/C},
$$
with probability 1, where $C>0$ is a certain (explicitly determined)
constant depending on $z,z'$. Roughly speaking, this means that the
Thoma parameters decay with the rate of a geometric progression. In
section 6, we produce the Sturm--Liouville differential operators
which commute with the integral operators given by the Whittaker
kernel and the stationary kernels of section 4. In section 7, we
compare the processes $\PPz^+$ with the Poisson--Dirichlet processes
$\PDt$ \cite{Ki}. 

\subhead 1. The fermion point processes \endsubhead
Let $\PP$ be a point process on a space $X$ and
$\rho_n(x_1,\dots,x_n)$ denote its correlation functions relative to a
reference measure $\mu$ on $X$, $n=1,2,\dots$. By a test set $A\subset
X$ we mean a Borel set such that the first correlation function
$\rho_1$ is integrable on $A$; this means that $A$ intersects the
random configuration at a finite number of points with probability 1.

In many concrete situations, it turns out that the functions $\rho_n$
are given by a determinantal formula,  
$$
\rho_n(x_1,\dots,x_n)=\det[K(x_a,x_b)]_{a,b=1}^n\,,
$$
where $K(x,y)$ is a kernel on $X\times X$ not depending on $n$.

The processes with determinantal correlation functions seem to be
important enough to deserve a special name. Though key examples
of such processes were already considered in the 1st edition (1967)
of Mehta's book on random matrices \cite{Me} and some earlier papers,
the first (to our knowledge) general discussion appeared in Macchi's
paper \cite{Ma1}. In her works and in the book \cite{DVJ} these
processes are called the {\it fermion point processes\/}, and we
shall adopt this terminology. 

Let us list some general properties of the fermion processes (see
\cite{Ma1, Ma2} and \cite{DVJ}, Example 5.4(c) and Exercises 5.4.7--9).

If the reference measure $\mu$ is replaced by an equivalent one,
$\mu\mapsto f\mu$, where $f$ is a strictly positive function, then
the kernel must be divided by $\sqrt{f(x)f(y)}$.  

Note that a determinantal correlation function vanishes when some of
the arguments coincide; this means that the points of the random
configuration have some repulsion properties.

Given a kernel $K(x,y)$, the following conditions ensure the
existence of a fermion process, see \cite{Ma1}:

\roster
\item"{(*)}" The functions $\rho_n$ as defined above are nonnegative.
\item"{(**)}" For a test set $A\subset X$, let $K_A$ denote the
integral operator whose kernel $K_A(x,y)$ is defined as the
restriction of $K(x,y)$ to $A\times A$. It is required that the norm
of $K_A$ in $L^2(A,\mu)$ be strictly less than 1 for any test set $A$.
\endroster

In many cases the kernel $K(x,y)$ turns out to be symmetric, so that
$K_A$ is a nonnegative self--adjoint operator in the Hilbert space
$L^2(A,\mu)$. But there also exist interesting examples of non
symmetric kernels (see, e.g., \cite{B2}). Even for symmetric kernels,
a direct verification of the above sufficient conditions can be
difficult. \footnote{It greatly simplifies when $K(x,y)$ is a
translation invariant kernel, say, on the real axis, see section 4 
below.} But $K(x,y)$ often arises as a limit of kernels which
certainly satisfy $(*)$ and $(**)$. In such a situation we can
conclude at least that $K(x,y)$ satisfies $(*)$ and $\Vert
K_A\Vert\le1$, which is a weak form of (**). 

Let us denote by $\pi^{(A)}_n(x_1,\dots,x_n)$ the finite--dimensional
distribution functions of $\PP$. Here $n=1,2,\dots$, $A\subset X$ is
a test set and, by definition,
$$
\gather
\pi^{(A)}_n(x_1,\dots,x_n)\mu(dx_1)\dots \mu(dx_n)
=\Prob\{\text{exactly $n$ points in $A$,} \\
\text{one point located in each of
the infinitesimal regions $dx_i$}\}.
\endgather
$$
By the well--known inclusion--exclusion principle, the $\pi$ functions
can be expressed through the $\rho$ functions as follows
$$
\pi^{(A)}_n(x_1,\dots,x_n)=\sum_{k=0}^\infty
\frac{(-1)^k}{k!}\int_{A^k}\rho_{n+k}(x_1,\dots,x_n,y_1,\dots,y_k)
\mu(dy_1)\dots\mu(dy_k),
$$
see \cite{DVJ}, section 5.4. For the fermion processes this relation
takes the following form
$$
\pi^{(A)}_n(x_1,\dots,x_n)=\Det(1-K_A)\det[L_A(x_a,x_b)],
$$
where $K_A$ is the integral operator whose kernel $K_A(x,y)$ is
obtained by restricting the kernel $K(x,y)$ to $A\times A$,
$\Det$ is the Fredholm determinant $\Det(1-\la K_A)$ evaluated at
$\la=1$, and $L_A(x,y)$ is the kernel of the operator
$L_A:=K_A(1-K_A)^{-1}$. In particular, the probability $\pi^{(A)}_0$
that $A$ is empty of points of the random configuration equals
$\Det(1-K_A)$.  

For a fermion process on the line, there is a relation between the
Fredholm determinant and the probability distribution of the
spacings, see, e.g., \cite{TW4}, first formula after \thetag{5.35}.
For other relations involving the Fredholm determinant, see \cite{Me,
TW1}. 

Now let us give examples of the kernels $K(x,y)$ originated from
concrete problems.

The most known is the {\it sine kernel}
$$
K(x,y)=\frac{\sin\pi(x-y)}{\pi(x-y)}, \qquad x,y\in\R.
$$
It appears in scaling limit of various random matrix ensembles
``in the bulk of the spectrum'', see \cite{Me, NW1, TW1}. Since the
kernel is translation invariant, the corresponding process is
stationary. The sine kernel can be included into a more general
family of translation invariant kernel, see section 4 below.
\footnote{About the meaning of the factor $\pi$, see the first
comment to Proposition 4.2.}

Taking scaling limit ``at the edge of the spectrum'' leads to
other kernels: the {\it Airy kernel}
$$
K(x,y)=\frac{Ai(x)Ai'(y)-Ai(y)Ai'(x)}{x-y}\,,
\qquad x,y\in\R
$$
(where $Ai(\cdot)$ is the Airy function) and the {\it Bessel kernel}
$$
K(x,y)=\frac{J_\al(x^{\frac12})y^{\frac12}J'_\al(y^{\frac12})
-J_\al(y^{\frac12})x^{\frac12}J'_\al(x^{\frac12})}{2(x-y)}\,,
\qquad x, y>0
$$
(where $J_\al(\cdot)$ is the Bessel function of order $\al>-1$), see
\cite{F, NW2, TW2, TW3}. 

The results of Part II, section 3, lead to one more kernel expressed
through special functions, the {\it Whittaker kernel}
$$
K(x,y)=\const\cdot(xy)^{-\frac12}\,
\frac{W_{\kappa,\mu}(x)W_{\kappa-1,\mu}(y)
-W_{\kappa,\mu}(y)W_{\kappa-1,\mu}(x)}{x-y}\,,
\qquad x,y>0
$$
(where $W_{\kappa,\mu}(\cdot)$ is the Whittaker function, a version
of the confluent hypergeometric function, see \cite{E1}). See also
section 2. We did not encounter this kernel in literature, maybe
this is a new example. 

Note that for the process with the sine kernel the points are
accumulated near $\pm\infty$; for the Airy kernel --- near $-\infty$;
for the Bessel kernel --- near $+\infty$; and for the Whittaker
kernel --- near 0.

An important problem is studying the Fredholm determinant
$\Det(1-\la K_A)$ (the test set $A$ being an interval or a finite union
of intervals), in particular, its asymptotics as an end of an
interval tends to an accumulation point. See \cite{TW1--5}.

\subhead 2. Lifting and the Whittaker kernel \endsubhead
Given a point process $Q$ on $(0,1]$ and a probability distribution $\si$
on $(0,+\infty)$, we can construct a new point process $\wtQ$, living
on $(0,+\infty)$, as follows. We take the random configuration $\xi$
corresponding to $Q$ and multiply it by an independent scale factor
$s$ distributed according to $\si$. Choose as $\si$ the gamma
distribution with density 
$$
\frac{\si(ds)}{ds}=\frac{s^{t-1}e^{-s}}{\Ga(t)}\,, \qquad t>0.
$$
Following Part II, section 3, we shall call $\wtQ$ the {\it
lifting\/} of $Q$ with parameter $t$.

Let us denote by $\PPz^+$ the restriction of the process $\PPz$ to
$(0,1]\subset I$. The process $\PPz^+$ governs the random behavior of the
Thoma parameters $\al_1,\al_2,\dots$ and neglects the parameters
$\be_1,\be_2,\dots$. On the contrary, to focus on the beta part of
the Thoma parameters it suffices to replace $z,z'$ by $-z,-z'$. 

In Part II, Theorem 3.3.4, we proved the following
result: 

\proclaim {Theorem 2.1} Application of lifting with parameter $t=zz'$
to the process $Q=\PPz^+$ gives a fermion process on $(0,+\infty)$.
Its kernel $K(x,y)$ is the Whittaker kernel
$$
\frac{(xy)^{-1/2}}{\Ga(z)\Ga(z')}\cdot
\frac{W_{\frac{z+z'+1}2,\frac{z-z'}2}(x)
W_{\frac{z+z'-1}2,\frac{z-z'}2}(y)
-W_{\frac{z+z'+1}2,\frac{z-z'}2}(y)
W_{\frac{z-z'-1}2,\frac{z-z'}2}(x)}{x-y}\,.
$$
\endproclaim

See section 6.1 in \cite{E1} for the definition of the Whittaker
function $W_{\kappa,\mu}$. 

\example{Comments} 1) Comparison with Theorems 2.2.1 and 2.4.1 from Part
II shows that lifting greatly simplifies the structure of the
expressions for the correlation functions. For instance, the
dimension of integrals involved in the description of the unlifted
correlation functions $\rz_n$ grows with $n$ while for the lifted
correlation functions $\widetilde\rho^{(zz')}_n$ we need only a
one--dimensional integral implicit in the definition of the Whittaker
function. 

2) By Proposition 3.1.1 (Part II), the passage from $\rz$ to
$\widetilde\rho^{(zz')}_n$ is effectued by the integral transform
$\Cal L^t$ defined in Proposition 3.2.1. This transform is
readily reduced to one--dimensional Laplace transform and so, in
principle, can be inverted via the Laplace inversion formula. This
implies that the lifted process retains the whole information about
the initial process.  

3) Let $(\al_1>\al_2>\dots)$ be the random configuration of the process
$\PPz^+$ and $(\widetilde\al_1>\widetilde\al_2>\dots)$ be the similar
object for the lifted process. The distribution functions of $\al_1$
and $\widetilde\al_1$ are also related by the transform $\Cal L^t$
and the same is true for joint distributions of any finite number of
coordinates.  

4) Note that the distribution of $\widetilde\al_1$ is
given by the Fredholm determinant: 
$$
\Prob\{\widetilde\al_1<\tau\}=\Det(1-K_{(\tau,+\infty)}),\qquad
\tau>0, 
$$
where $K_{(\tau,+\infty)}(x,y)$ is the restriction of the Whittaker
kernel to $(\tau,+\infty)$.

5) The transform $\Cal L^t$ has a simple meaning in the language of
moments. Its application, say, to a one-dimensional distribution
results in multiplying the $m$th moment by $(t)_m$,
$m=1,2,\dots$. Thus, if we calculate numerically a few moments of the
random variable $\widetilde\al_1$ then we immediately get the
corresponding moments of $\al_1$. 

6) The fact that application of lifting can simplify a point process is
also demonstrated on the example of the Poisson--Dirichlet process
whose lifting is simply a Poisson process (Part II, Proposition
3.1.2). Note that from this result one can very easily get Griffiths'
formulas \cite{G1, G2} for the mean values of the random coordinates
$x_1>x_2>\dots$ distributed according to the Poisson--Dirichlet law. 
\endexample

As a simple application of Theorem 2.1 we shall prove the following
result. 

\proclaim{Proposition 2.2} Consider the probability space $(\Om,\Pz)$
where $\Om$ is the Thoma simplex and $\Pz$ is one of the measures of
the principal or complementary series {\rm(}see Part I\/{\rm)}. Let
us view the Thoma parameters as random variables defined on this
probability space and let the symbol $\Bbb E$ mean expectation. Then
we have 
$$
\Bbb E\left(\sum_{i=1}^\infty\alpha_i\right)=
\frac{\sin\pi z\cdot\sin\pi z'}{\pi\sin\pi(z-z')}
\left[(z-z')\frac{z+z'-1}{2zz'}+\psi(-z')-\psi(-z)
\right]\,,\tag2.1
$$
where $\psi(a)=\Gamma'(a)/\Gamma(a)$.

Similarly,
$$
\Bbb E\left(\sum_{i=1}^\infty\beta_i\right)=
\frac{\sin\pi z\cdot\sin\pi z'}{\pi\sin\pi(z-z')}
\left[(z'-z)\frac{z+z'+1}{2zz'}+\psi(-z)-\psi(-z')
\right]\,.\tag2.2
$$
\endproclaim

\demo{Proof} By Proposition 4.6 of Part I, the symmetry map
$(\al,\be)\mapsto(\be,\al)$ of the Thoma simplex takes the measure
$\Pz$ to the measure $P_{-z,-z'}$. On the other hand, the right--hand
sides of the formulas \tht{2.1}, \tht{2.2} differ exactly by change
of sign in $z,z'$. So, the both formulas are equivalent, and it
suffices to check one of them, say, \tht{2.1}.

Let $\rho_1(x)$ denote the density function of the process $\PPz$. By
the very definition of the density function (see Part I, \S4),
$$
\Bbb E\left(\sum_{i=1}^\infty\al_i\right)=
\int_0^1 x\rho_1(x)dx. \tag2.3
$$
Since we know various expressions for the density function (see
Part I, Theorem 5.10, Theorem 5.12; Part II, Corollary 2.4.2),
we could try to employ one of them to calculate the integral \tht{2.3}
explicitly. However, this does not seem to be easy, so we have
preferred to use a roundabout way -- reduction to the lifted process.

Let $\widetilde\al_1, \widetilde\al_2,\dots$ stand for the
``lifted'' random variables $\al_1,\al_2,\dots$ and
$\widetilde\rho_1$ be the density function of the 
lifted process (see section 3.1 in Part II). Similarly to \tht{2.3}
we have 
$$
\Bbb E\left(\sum_{i=1}^\infty\widetilde\al_i\right)=
\int_0^\infty x\widetilde\rho_1(x)dx. \tag2.4
$$
On the other hand, it follows from the definition of lifting (see
also Comment 5 above) that
$$
\int_0^1 x\rho_1(x)dx=\frac1t\int_0^\infty x\widetilde\rho_1(x)dx,
$$
where, as usual, $t=zz'$. 

Let $K(x,y)$ be the Whittaker kernel (Theorem 2.1). Then
$$
\widetilde\rho_1(x)=K(x,x),\qquad x>0.
$$
Let us abbreviate 
$$
\ka=\frac{z+z'+1}{2}, \quad \mu=\frac{z-z'}2,\tag2.5
$$
and assume $x>0$. Then, applying the l'Hospital rule, we get 
$$
x\widetilde\rho_1(x)=
\frac1{\Ga(z)\Ga(z')}
(W'_{\ka,\mu}(x)W_{\ka-1,\mu}(x)-W_{\ka,\mu}(x)W'_{\ka-1,\mu}(x)). \tag2.6
$$
Let us employ the relation
$$
W'_{\kappa,\mu}=\left(-\frac12+\frac\kappa{x}\right)W_{\kappa,\mu}+
\frac{(\frac12-\kappa+\mu)(\frac12-\kappa-\mu)}{x}W_{\kappa-1,\mu}\,,
$$
which follows from formulas 6.9 (2) and 6.6 (1) in \cite{E1}. Then we
get from \tht{2.5} and \tht{2.6} 
$$
\multline
\frac1t x\wtrho_1(x)\\
=\frac{W_{\ka,\mu}(x) W_{\ka-1,\mu}(x)+
zz' W^2_{\ka-1,\mu}(x)-(z-1)(z'-1)W_{\ka,\mu}(x) W_{\ka-2,\mu}(x)}
{zz'\cdot\Gamma(z)\Gamma(z')\cdot x}\,.
\endmultline
$$

The integral of the above expression over $(0,+\infty)$ can be found
by making use of the table integrals
$$
\multline
\int_0^\infty x^{-1}W_{\ka_1,\mu}(x)W_{\ka_2,\mu}(x)dx
=\frac\pi{(\ka_1-\ka_2)\sin(2\pi\mu)}\\
\times\left[
\frac1{\Ga(\frac12-\ka_1+\mu)\Ga(\frac12-\ka_2-\mu)}
-\frac1{\Ga(\frac12-\ka_1-\mu)\Ga(\frac12-\ka_2+\mu)}
\right]
\endmultline
$$
and
$$
\int_0^\infty x^{-1}W^2_{\ka,\mu}(x)dx=
\frac\pi{\sin2\pi\mu}\cdot
\frac{\psi(\frac12-\ka+\mu)-\psi(\frac12-\ka-\mu)}
{\Ga(\frac12-\ka+\mu)\Ga(\frac12-\ka-\mu)}\,,
$$
which can be found in \cite{PBM}, section 2.19.23, formulas 3 and 4.
Note that these two integrals are convergent provided that $|\Re\mu|<\frac12$,
which indeed holds in our case, because of the restrictions imposed
on the parameters $z,z'$, see \S2 in Part I.

Finally, we also need the relation 
$$
\psi(a+1)=\psi(a)+\frac1a\,,
$$
see \cite{E1}, 1.7 (8).

Then, after elementary calculations we get the desired formula
\tht{2.1}. \qed
\enddemo

\example{Remark 2.3} Assume that $z'=-z$; according to the assumptions on
the parameters $z,z'$, this means that $z=-z'$ is purely imaginary.
Then the measure $\Pz$ is stable under the symmetry map transposing
the $\al$'s and the $\be$'s, so that the expressions \tht{2.1} and
\tht{2.2} must be equal to $1/2$. One can check that this is indeed
 the case by making use of the relation 
$$
\psi(a)-\psi(-a)=-\pi\operatorname{ctg}(\pi a)-\frac1a\,,
$$
which can be found in \cite{E1}, 1.7 (11).

Further, we know that $\sum(\al_i+\be_i)=1$ almost surely with
respect to $\Pz$ (Part I, Theorem 6.1). It follows that the sum of
the expressions \tht{2.1} and \tht{2.2} must be identically equal to 1.
Again, this can be readily verified by making use of the above
relation. 
\endexample

\example{Remark 2.4} Note that for certain special values of the
parameters $z,z'$ the Whittaker kernel degenerates to the
Christoffel--Darboux kernel for the Laguerre polynomials.

Specifically, let $L^{2\mu}_N$ stand for the $N$th Laguerre
polynomial with the weight function $x^{2\mu}e^{-x}$ on $\R_+$, where
$2\mu>-1$; the normalization is that of \cite{E1}. We have
$$
x^{-\frac12}W_{\mu+N+\frac12, \mu}(x)=
(-1)^N N! x^\mu e^{-\frac x2}L^{2\mu}_N(x),
$$
see, e.g. \cite{E1}, 6.9 (36).

Assume $N-1<z,z'<N$, where $N=1,2,\dots$, and let $z'$ tends to $N$
while $z$ remains fixed. Denote $2\mu=z-N$; then $-1<2\mu<0$. Using
the above formula we readily get that the limit of the Whittaker
kernel as $z'\to N$ is equal to
$$
\gather
\frac{N!}{\Ga(N+2\mu)}\,(xy)^\mu e^{-\frac{x+y}2}\,\frac
{L^{2\mu}_{N-1}(x)L^{2\mu}_N(y)-L^{2\mu}_{N-1}(y)L^{2\mu}_N(x)}{x-y}\\
=(xy)^\mu e^{-\frac{x+y}2}\,\sum_{i=0}^{N-1}
\frac{L^{2\mu}_i(x)L^{2\mu}_i(y)}
{\int (L^{2\mu}_i(x))^2x^{2\mu}e^{-x}dx}\,.
\endgather
$$
The latter expression coincides with the kernel of the projection in
the Hilbert space $L^2(R_+, \,dx)$ on the linear span of the functions
$x^{i+\mu}e^{-\frac x2}$, where $i=0,\dots,N-1$; this is
exactly the kernel associated with the ``$N$-point Laguerre
polynomial ensemble'', see  [FK], [Br], [NW1].  

Finally, note that the restriction $\mu<0$, which comes from the
assumption $N-1<z,z'<N$, is inessential, because there exists a
natural ``degenerate series'' of the measures $\Pz$ with the
parameters $z'=N$, $z>N-1$. 
\endexample

\subhead 3. The tail process \endsubhead
Let $Q$ be a point process on $(0,1]$ or on $(0,+\infty)$ and
$\rho_1(x)$ be its first correlation function. We assume that
$\rho_1$ is integrable on the right of any $\ep>0$ and
nonintegrable on $(0,\ep)$, so that the points are accumulated to
0. Consider the mapping
$$
(0,1]\,\to\, [0,+\infty),\qquad
x\mapsto\xi:=\int_x^1\rho_1(y)dy,
$$
and let $\whQ$ be the image of the process $Q$ (or rather of its
restriction to $(0,1]$, in case $Q$ is defined on the whole ray)
under that mapping. Then $\whQ$ is a point process on 
$[0,+\infty)$ and its first correlation measure coincides with
Lebesgue measure. 

Further, for any $\tau\ge0$, let $\whQ_\tau$ be the process on
$[-\tau,+\infty)$ obtained from $\whQ$ by the shift
$\xi\mapsto\xi-\tau$. We let $\tau\to+\infty$ and assume that there
exists a point process $\whQ_\infty$ on the whole axis $\R$ such that
the limit
$$
\lim_{\tau\to+\infty}\whQ_\tau=\whQ_\infty
$$ exists in a reasonable sense. Then we shall say that $\whQ_\infty$ is
the {\it tail process\/} for $Q$.

Of course, the exact meaning of the limit above has to be precised.
We shall be content with the following type of
convergence: for any $n$, the $n$th correlation function of
$\whQ_\tau$ tends, as $\tau\to+\infty$, to the $n$th correlation
function of $\whQ_\infty$, uniformly on compact sets in $\R$.
Perhaps, the definition can be elaborated. But anyway, the idea is
clear: we restrict the initial process  to a
small interval $(0,\ep)$, next rescale it to make the density function
constant, and then look at the asymptotics as $\ep\to0$.

In the examples considered below the correlation functions of $Q$ can
be represented in the form
$$
\rho_n(x_1,\dots,x_n)=\frac{C^n}{x_1\dots x_n}(f_n(x_1,\dots,x_n)+o(1)),
\qquad x_1,\dots,x_n>0,
$$
where $C>0$ is a constant not depending on $n$, $f_n(x_1,\dots,x_n)$
is a continuous homogeneous function,
$$
f_n(rx_1,\dots,rx_n)=f_n(x_1,\dots,x_n)\qquad r>0,
$$
and the rest term, denoted as $o(1)$, tends to zero as
$\max\{x_1,\dots,x_n\}\to0$. Note that the function $f_1$ should be a
constant, and we choose $C$ in such a way that $f_1(\cdot)\equiv1$.

In such a situation, we make a change of variables $x_i\mapsto\xi_i$,
where $x_i=e^{-\xi_i/C}$. In the new variables, the correlation functions
take the form
$$
\rho'_n(\xi_1,\dots,\xi_n)=g_n(\xi_1,\dots,\xi_n)+o(1),
\qquad \xi_1,\dots,\xi_n\in\R,
$$
where the rest term $o(1)$  tends to zero
as $\min\{\xi_1,\dots,\xi_n\}\to+\infty$ and
$$
g_n(\xi_1,\dots,\xi_n):=f_n(e^{-\xi_1/C},\dots,e^{-\xi_n/C})
$$
is a translation invariant function. Consequently, the desired
``tail'' correlation functions have the form
$$
\widehat\rho(\xi_1,\dots,\xi_n)=g_n(\xi_1,\dots,\xi_n).
$$

As illustration, examine first the Poisson--Dirichlet process.

\proclaim{Proposition 3.1} The tail process for the
Poisson--Dirichlet process $\PDt$ is the  Poisson process on
$\R$ with constant density 1.
\endproclaim

\demo{Proof} Recall that the correlation functions of $\PDt$ are
given by Watterson's formula
$$
\rho_n(x_1,\dots,x_n)=\frac{t^n(1-x_1-\dots-x_n)_+^{t-1}}
{x_1\dots x_n}\,,
$$
see \cite{W} and Part I, Corollary 7.4. In particular, the first
correlation function is 
$$
\rho_1(x)=\frac{t(1-x)^{t-1}}x\,.
$$

These correlation functions fit into the above scheme with $C=t$
and all the functions $f_n$ identically equal to 1. 
It follows that the ``tail'' correlation functions $\widehat\rho_n$ are
identically equal to 1, which corresponds to the standard Poisson
process. \qed 
\enddemo

As in section 2 above, let $\PPz^+$ denote the restriction of the
process $\PPz$ to $(0,1]\subset I$.

\proclaim{Theorem 3.2} Take as $Q$ the process $\PPz^+$ or its
lifting with parameter $t=zz'$. In both cases the tail process
$\wtQ$ is a fermion process on $(0,+\infty)$ with the same translation
invariant kernel $K(\xi,\eta)$, which has the following form.  

$\bullet$ For the principal series, when $z'=\bar z$ and $z$ is not real,
$$
K(\xi,\eta)=\frac{B\sin A(\xi-\eta)}{A\sh B(\xi-\eta)}\,,
$$
where
$$
B=\frac{\pi\sin\pi(z-z')}{2(z-z')\sin\pi z\cdot\sin\pi z'}\,>0\,,
\qquad A=\pm i(z-z')B.
$$

$\bullet$ For the supplementary series, when $m<z,z'<m+1$ for a certain
$m\in\Z$ and $z\ne z'$, 
$$
K(\xi,\eta)=\frac{B\sh A(\xi-\eta)}{A\sh B(\xi-\eta)}\,,
$$
where $B$ is given by the same formula and $A=\pm(z-z')B$.

$\bullet$ Finally, on the intersection of the both series, when
$z=z'\in\R\setminus\Z$, the kernel is given by the limit expression
$$
K(\xi,\eta)=\frac{B(\xi-\eta)}{\sh B(\xi-\eta)}\,,
$$
where
$$
B=\frac{\pi^2}{2\sin^2\pi z}\,.
$$
\endproclaim

\demo{Proof} For the lifted process, the behavior of the
correlation functions near zero is given by the asymptotics of the
Whittaker kernel as described in Part II, Theorem 4.1.1. For the
process $\PPz^+$ itself this requires knowledge of the asymptotics of the
multivariate Lauricella functions of type $B$; the final result is
described in Part II, Theorem 4.3.1. According to these theorems, the
correlation functions, both in lifted and non lifted case, fit into
the above scheme with the same constant $C$ and the same functions
$f_n$. 

Specifically, we have 
$$
C=\cases\dfrac{(z-z')\sin\pi z\cdot\sin\pi z'}
{\pi\sin\pi(z-z')},\quad z'\ne z\\
\dfrac{\sin^2\pi z}{\pi^2}, \quad z'=z\in\R\setminus\Z
\endcases
$$  
and
$$
f_n(x_1,\dots,x_n)=\det[K'(x_i,x_j)],
$$
where
$$
K'(x,y)=\cases\dfrac1{z-z'}\cdot
\dfrac{(x/y)^{\frac{z-z'}2}-(x/y)^{\frac{z'-z}2}}
{(x/y)^{\frac12}-(x/y)^{-\frac{1}2}}\,,\quad z'\ne z\\
\dfrac{\ln x-\ln y}{(x/y)^\frac12-(x/y)^{-\frac12}}\,, 
\quad z=z'\in\R\setminus\Z.
\endcases
$$

Passing from the functions $f_n(x_1,\dots,x_n)$ to the functions
$g_n(\xi_1,\dots,\xi_n)$ as described above we get the expressions
indicated in the statement of the theorem. \qed
\enddemo

\example{Remark 3.3} Let us notice two remarkable features of the
expression for the kernel in Theorem 3.2. 

$\bullet$ First, the kernel does not change under the
transform $(z,z')\mapsto(-z,-z')$. This implies that the
tail properties of the $\al_j$'s are  the same as that of the $\be_j$'s. 

$\bullet$ Second, the kernel does not change  under the shift 
$(z,z')\mapsto(z+1,z'+1)$. This  means a quite surprising periodicity
of the tail process with respect to the parameters $z,z'$. One can ask
whether this phenomenon is somehow related to degeneration of $\Pz$
at integer values \cite{KOV}. 
\endexample

\subhead 4. The $\sin/\sh$ and $\sh/\sh$ kernels \endsubhead
Here we shall examine in more detail the kernels that appeared in Theorem
3.2.  These are the stationary kernels of the form
$$
K(x,y)=\frac{B\sh A(x-y)}{A\sh B(x-y)}
$$
where $B$ is real and $A$ is either real or pure imaginary (here and
below we use the letters $x,y$ instead of $\xi,\eta$). There are
two main types and two limit types: 

1) The ``$\sin/\sh$ kernel'', 
$$
K(x,y)=\frac{B\sin A(x-y)}{A\sh B(x-y)}, 
\qquad B>0, \quad A\in\R, \quad A\ne0.
$$

2) The ``$\sh/\sh$ kernel'', 
$$
K(x,y)=\frac{B\sh A(x-y)}{A\sh B(x-y)}, 
\qquad B>0, \quad A\in\R, \quad A\ne0.
$$

3) The limit case $B=0$:
$$
K(x,y)=\frac{\sh A(x-y)}{A(x-y)}, 
\qquad  A\in\R, \quad A\ne0.
$$

4) The limit case $A=0$:
$$
K(x,y)=\frac{B(x-y)}{\sh B(x-y)}, \qquad B>0.
$$

In all the cases we normalized the kernels so that $K(x,x)\equiv1$; by a
change of variable, $x\mapsto\operatorname{const}\cdot x$, we could
replace 1 by an arbitrary constant.

\proclaim{Proposition 4.1} Let $K(x,y)=k(x-y)$ be a translation invariant
kernel on $(\R,dx)$ and assume that $k(\cdot)$ is
the inverse Fourier transform of an integrable function $\whk(\cdot)$,
$$
k(x)=\frac1{2\pi}\int_{-\infty}^{+\infty}e^{-ixy}\whk(y)dy,
$$
such that $0\le\whk(y)\le1$ for all $y\in\R$ and $\whk(y)<1$ when
$|y|$ is large enough. 

Then the conditions $(*)$ and $(**)$ stated in section 1 are
satisfied, so that the kernel corresponds to a fermion process on
$\R$. 
\endproclaim

\demo{Proof} Since $\whk$ is nonnegative,  the function $k$ is
Hermitian--symmetric and nonnegative definite, whence the condition
(*) is satisfied. 

Let $K$ denote the integral operator in $L^2(\R,dx)$ with the kernel
$K(x,y)$. The image of $K$ under the Fourier transform is the operator of
multiplication by the function $\whk$. This implies that $0\le K\le1$. 

It remains to check that $K_A$ is strictly less than 1 for any
test set $A\subset\R$, i.e., for any bounded $A$; without loss of
generality one can assume that $A$ is a finite interval $[a,b]$. The
function $k$ being continuous, the operator $K_A$ is a compact
Hermitian operator. Assume $K_A\phi=\phi$ for a function $\phi\in
L^2(\R)$. Then $\phi$ is concentrated on $A$ and $K\phi=\phi$. Taking
the Fourier transform we see that the Fourier image $\widehat\phi$ is
concentrated on the region where $\whk(\cdot)=1$. By the assumption,
this region is bounded, so that both $\phi$ and $\widehat\phi$ are
compactly supported, which implies $\phi\equiv0$. This means that
$K_A<1$, which completes the proof. \qed
\enddemo

Let $k(x)=B\sin Ax/A\sh Bx$ or $k(x)=B\sh Ax/A\sh Bx$, where in the
latter case we assume $|A|<B$ (otherwise $k(x)$ certainly does not have
the required form). Then the Fourier transform $\whk(y)$ is given by
the formula
$$
\whk(y)=\frac{\pi\sh(\pi A/B)}{A[\ch(\pi A/B)+\ch(\pi y/B)]}
$$
or
$$
\whk(y)=\frac{\pi\sin(\pi A/B)}{A[\cos(\pi A/B)+\ch(\pi y/B)]}\,,
$$
respectively. These formulas (which are related by analytic
continuation with respect to the parameter $A$) can be found, e.g.,
in  \cite{E2}, section 1.9, formula 14. 

\proclaim{Proposition 4.2} The above four stationary kernels
generate a fermion point process if the parameters $A,B$
satisfy the following conditions, respectively.

{\rm1)} The ``$\sin/\sh$ kernel'':
$$
\th\frac{\pi|A|}{2B}\le\frac{|A|}\pi\,.
$$

{\rm2)} The ``$\sh/\sh$ kernel'':
$$
0\le\,\tg\frac{\pi|A|}{2B}\le\frac{|A|}\pi\,,
\qquad \frac{|A|}B\,<\,1.
$$

{\rm3)} The limit case $B=0$: 
$$
A\ge\pi.
$$

{\rm4)} The limit case $A=0$: 
$$
B\ge\pi^2/2.
$$
\endproclaim

\demo{Proof} The above expressions for $\whk(y)$ correspond to the
first two cases. For the remaining two cases we get $\whk(y)$ by an
obvious limit transition. We must find conditions on $A,B$ under
which $\whk(y)$ satisfies the two inequalities of Proposition 4.1.
This is done by an elementary calculation. \qed
\enddemo

\example{Comments} 1) In our scheme, the famous sine kernel
corresponds to the degenerate case $B=0$ and the minimal allowed value
$A=\pi$.  

2) The ``$\sin/\sh$'' kernel with $A=\pi$ and
arbitrary $B$ appeared in the papers \cite{BCM, MCIN}.

3) For the principal series, when $z=a+ib$, $z'=a-ib$, where $a,b$ are
real and $b$ is nonzero, the kernel of the tail process is of type
``$\sin/\sh$'' with the parameters
$$
A=\frac{\pi\sh(2\pi b)}{\ch(2\pi b)-\cos(2\pi a)}\,,\quad
B=\frac{\pi\sh(2\pi b)}{2b(\ch(2\pi b)-\cos(2\pi a))}\,.
$$
The inequality imposed on $A,B$ becomes evident in terms of $a,b$: 
$$
\frac{\ch(2\pi b)-\cos(2\pi a)}{\ch(2\pi b)+1}\,\le1.
$$
Note that we do not get all the allowed couples $(A,B)$.

4) For the complementary series, when $z$ and $z'$ are real such that
$m<z,z'<m+1$ for a certain $m\in\Z$ and $z\ne z'$, the tail kernel is
of type ``$\sh/\sh$'' with the parameters
$$
A=\frac{\pi\sin\pi(z-z')}{2\sin\pi z\cdot\sin\pi z'}\,, \quad
B=\frac{\pi\sin\pi(z-z')}{2(z-z')\sin\pi z\cdot\sin\pi z'}\,.
$$
The first inequality imposed on $A,B$ turns into the evident one:
$$
0\le\,\frac{2\sin\pi z\cdot\sin\pi z'}{\cos\pi(z-z')+1}=
\frac{\cos\pi(z-z')-\cos\pi(z+z')}{\cos\pi(z-z')+1}\,\le1.
$$
The second inequality is also evident, because
$$
\frac{|A|}B=|z-z'|<1.
$$
Again, we do not get all the allowed couples $(A,B)$. 

5) For the intersection of the both series, when $z=z'\in\R\setminus\Z$,
the kernel is of limit type with the parameters
$$
A=0, \quad B=\frac{\pi^2}{2\sin^2\pi z}\,,
$$
and the inequality on $B$ takes the form
$$
\frac{\pi^2}{2\sin^2\pi z}\,\ge\,\frac{\pi^2}2\,.
$$
As $\sin^2\pi z$ ranges over $(0,1]$, we get here all the allowed 
values of the parameter $B$. 

6) Consider the principal series ($z=a+ib$, $z'=a-ib$) and let
$|b|\to\infty$; then, irrespective to the behavior of $a$, the
corresponding tail kernel tends to the sine kernel
$\sin\pi(\xi-\eta)/\pi(\xi-\eta)$. 

\endexample

\subhead 5. Rate of decay of the Thoma parameters \endsubhead
\proclaim{Theorem 5.1} Let $\Pz$ be an arbitrary measure of principal
or complementary series. Then, for almost every point
$\om=(\al,\be)\in\Om$, with respect to the measure $\Pz$, there exist
the limits
$$
\lim_{j\to\infty}(\al_j)^{1/j} =
\lim_{j\to\infty}(\be_j)^{1/j}=e^{-C^{-1}},
$$
where $C$ is the same as in the proof of Theorem 3.2, i.e.,
$$
C=\cases\dfrac{(z-z')\sin\pi z\cdot\sin\pi z'}
{\pi\sin\pi(z-z')},\quad z'\ne z\\
\dfrac{\sin^2\pi z}{\pi^2}, \quad z'=z\in\R\setminus\Z
\endcases
$$  
\endproclaim

\demo{Proof}{\it Step\/} 1 (change of a variable). We shall prove that the
first limit exists and equals $e^{-C^{-1}}$. Then this will also imply
the claim concerning the second limit, because, on the one hand,
$\al$ and $\be$ change places as $z$ and $z'$ are multiplied by $-1$,
and, on the other hand, this does not affect the value of $C$.

So, we shall examine the point process $\PPz^+$ on $(0,1]$. It is
convenient to pass from $(0,1]$ to $[0,+\infty)$ via the map
$x\mapsto\xi=-\ln x$. Let $x_1\ge x_2\ge \dots$ be the random
configuration of the process $\PPz^+$ and 
$\xi_1\le\xi_2\le\dots$ be its image under this map.
(Actually, we know that the inequalities are strict (see Theorem
2.5.1 of Part II) but here this is unessential.)

We have
$$
\lim_{j\to\infty}(x_j)^{1/j}=e^{-C^{-1}}\,\Leftrightarrow\,
\lim_{j\to\infty}\frac{\xi_j}j=C^{-1}
$$
(here we mean limits with probability 1).

{\it Step\/} 2 (reduction to $\Nt$). Let $\xi_1\le\xi_2\le\dots$ be the
random configuration for a point process on $[0,+\infty)$ and $\Nt$
denote the number of points in $[0,\tau]$, where $\tau>0$ is
arbitrary. Then the following equivalence holds ($C>0$ is a
constant): 
$$
\lim_{j\to\infty}\frac{\xi_j}j=C^{-1}\,\Leftrightarrow\,
\lim_{\tau\to+\infty}\frac\Nt\tau=C,
$$
limits with probability 1.

Actually, in this claim, randomness is unessential; we shall prove it
for a fixed (nonrandom) configuration.  

Assume that $\xi_j/j$ tends to $C^{-1}$ with a certain $C>0$. Then
$\xi_j\to+\infty$. Further, for any $\tau>0$, there exists a unique
$j$ such that $\xi_j\le\tau<\xi_{j+1}$. Then $\Nt=j$ and
$$
\frac j{\xi_{j+1}}\,<\,\frac\Nt\tau\,\le\frac j{\xi_j}\,.
$$
As $\tau\to+\infty$, we have $j\to+\infty$, whence, by the
assumption, the both bounds tend to $C$. Consequently, $\Nt/\tau$
tends to $C$, too.

Conversely, assume that $\Nt/\tau$ tends to a certain $C>0$. Then it
follows that $\xi_j\to+\infty$ as $j\to\infty$. Further, fix an
arbitrary $\ep>0$ and remark that
$$
N_{\xi_j-\ep}<j\le N_{\xi_j}\,.
$$
This implies
$$
\frac{N_{\xi_j-\ep}}{\xi_j}<\frac j{\xi_j}
\le\frac{N_{\xi_j}}{\xi_j}\,.
$$
As $j\to\infty$, we have $\xi_j\to+\infty$, whence, by the
assumption, the both bounds tend to $C$. It follows that $j/\xi_j$
tends to $C$, too, so that $\xi_j/j$ tends to $C^{-1}$.

{\it Step\/} 3 (use of correlation functions). Consider a point process on
$[0,+\infty)$. As above, by $\Nt$ we denote the number of points in
$[0,\tau]$. Let $\rho_1$ and $\rho_2$ be the first and the second
correlation functions.

\proclaim{Lemma 5.2} Assume that $\rho_1$ and $\rho_2$ satisfy the
following conditions, where $\ep$ and $C$ are certain constants
($0<\ep<1$ and $C>0$) and $\tau\to+\infty$: 
$$
\gather
\int_0^\tau\rho_1(\xi)d\xi=C\tau+O(\tau^{1-\ep})\\
\int_0^\tau\int_0^\tau\rho_2(\xi,\eta)d\xi d\eta
=C^2\tau^2+O(\tau^{2-\ep}).
\endgather
$$
Then $\Nt/\tau\to C$.
\endproclaim

\demo{Proof of the lemma} We shall adapt Kingman's argument in
\cite{Ki}, section 4.2. By the definition of the correlation measures,
$$
\gather
\int_0^\tau\rho_1(\xi)d\xi=\E(\Nt)\\
\int_0^\tau\int_0^\tau\rho_2(\xi,\eta)d\xi d\eta=\E(\Nt(\Nt-1)).
\endgather
$$
By the assumptions of the lemma, this implies
$$
\gather
\E(\Nt)=C\tau+O(\tau^{1-\ep})\\
\E(\Nt^2)=\E(\Nt(\Nt-1)+\Nt)=C^2\tau^2+O(\tau^{2-\ep}).
\endgather
$$
It follows
$$
\E\left(\left(\frac\Nt\tau-C\right)^2\right)=
\frac{\E(\Nt^2)}{\tau^2}-2C\frac{\E(\Nt)}\tau+C^2=O(\tau^{-\ep}).
$$
By the Chebyshev inequality, for any $\de>0$,
$$
\Prob\left(\left|\frac\Nt\tau-C\right|\ge\de\right)\,\le\,
\frac{\operatorname{const}}{\de^2\tau^\ep}\,,
$$
where the constant does not depend on $\tau$. Taking $\tau=k^{2/\ep}$,
where $k=1,2,\dots$, we conclude that the series
$$
\sum_{k=1}^\infty\Prob
\left(\left|\frac{N_{k^{2/\ep}}}{k^{2/\ep}}-C\right|\ge\de\right)
$$
converges for any $\de>0$. By the Borel--Cantelli lemma,
$$
\lim_{k\to\infty}\frac{N_{k^{2/\ep}}}{k^{2/\ep}}=C
$$
with probability 1.

Finally, for an arbitrary $\tau>0$, define a natural $k$ from the
relation 
$$
k^{2/\ep}\le\tau<(k+1)^{2/\ep}.
$$
Then
$$
N_{k^{2/\ep}}\le\Nt\le N_{(k+1)^{2/\ep}}\,.
$$
Since
$$
\frac{(k+1)^{2/\ep}}{k^{2/\ep}}\,\sim\,1, \qquad k\gg1,
$$
we get
$$
\lim_{\tau\to+\infty}\frac\Nt\tau=C
$$
with probability 1. \qed
\enddemo

{\it Step\/} 4 (reduction to the lifted process). Remark that
multiplication of a sequence $(x_1>x_2>\dots)$ by a positive factor
does not affect on the limit behavior of $x_j^{1/j}$. By the very
definition of the lifting, it follows that in the claim of the
theorem, we may replace our process by its lifting. The only purpose
of this reduction is that below we may employ Theorem 4.1.1 (Part II)
for the lifted process instead of the parallel Theorem 4.3.1 whose
proof is more difficult.

{\it Step\/} 5 (estimation of correlation functions). Consider the
lifting of the process $\PPz^+$ and then make change of a variable
$x\mapsto \xi=-\ln x$. Denote by $\rho_n(\xi_1,\dots,\xi_n)$ the
correlation functions of the resulting process. By the above
discussion it remains to check that $\rho_1$ and $\rho_2$ obey
the assumptions of Lemma 5.2.

Theorem 4.1.1 of Part II implies that our correlation functions can
be written in the form
$$
\rho_n(\xi_1,\dots,\xi_n)=C^n\det[K(\xi_i,\xi_j)]+r_n(\xi_1,\dots,\xi_n),
\qquad \xi_1,\dots,\xi_n>0,
$$
where the constant $C$ is the same as in the statement of the theorem,
$$
\gather
K(\xi,\eta)=k(\xi-\eta),\\
k(\zeta)=\cases
\dfrac{\sh\left(\frac{z-z'}2\,\zeta\right)}
{(z-z')\sh\left(\frac12\,\zeta\right)}\,,\quad z'\ne z\\
\dfrac{\zeta}{\sh\left(\frac12\zeta\right)}\,,
\quad z=z'\in\R\setminus\Z,
\endcases
\endgather
$$
and the rest term admits the estimate
$$
r_n(\xi_1,\dots,\xi_n)=O(e^{-\de\min\{\xi_1,\dots,\xi_n\}})
$$
with a certain $\de>0$ (this $\de$ is equal to 1 for the principal
series with $z'\ne z$; to $1-|z-z'|$ for the complementary series
with $z\ne z'$; and can be any number strictly less than 1, for the
intersection of the both series). 

Note that $k(0)=1$ and the function $k(\zeta)$ is an even square
integrable function on the whole real axis.

We have
$$
\int_0^\tau r_1(\xi)d\xi=O(1), \qquad
\int_0^\tau\int_0^\tau r_2(\xi,\eta)d\xi d\eta=O(\tau).
$$
Using this we get
$$
\gather
\int_0^\tau \rho_1(\xi)d\xi=C\int_0^\tau k(0)d\xi\,+\,O(1)
=C\tau\,+\,O(1),\\
\int_0^\tau\int_0^\tau\rho_2(\xi,\eta)d\xi d\eta=C^2
\int_0^\tau\int_0^\tau\vmatrix 1 & k(\xi-\eta)\\
k(\eta-\xi) & 1\endvmatrix d\xi d\eta\,+\,O(\tau)\\
=C^2\tau^2-C^2
\int_0^\tau\int_0^\tau k^2(\xi-\eta) d\xi d\eta\,+\,O(\tau).
\endgather
$$
Finally,
$$
\int_0^\tau\int_0^\tau k^2(\xi-\eta) d\xi d\eta
=2\int_0^\tau d\xi\int_0^\xi k^2(\zeta)d\zeta=O(\tau),
$$
because $k(\zeta)$ is square integrable.

Thus, we have verified the assumptions of Lemma 5.2, which completes
the proof. \qed
\enddemo

\example{Remark 5.3} The same argument can be applied to the
Poisson--Dirichlet process $\PDt$. Here the final step 5 is much
easier, because of a simpler structure of the correlation functions.
We get in this way that for $\PDt$, $(x_j)^{1/j}$ tends to
$e^{-t}$, the result originally obtained (for $t=1$) in \cite{VS} by a
quite different way.  We can also use lifting, as suggested on step
4, which provides a quick reduction to the law of large numbers for
the Poisson process, see \cite{Ki}, section 4.2.  
\endexample

Thus, both for $\PPz^+$ and $\PDt$, the rate of decay of the $x_j$'s
is of the same type. 

\example{Remark 5.4} Let $\Om'$ denote the subset of points 
$\om=(\al,\be)\in\Om$ such that the limit $F(\om):=\lim(\al_j)^{1/j}$
exists. Clearly, this is a Borel subset. According to Theorem 5.1,
$\Om'$ is of full measure with respect to any $\Pz$ and the function
$F$ takes constant values almost everywhere, with the constant
depending on $z,z'$ in a nontrivial way. This agrees with the fact that
the measure $\Pz$ are pairwise disjoint. However, this does not
provide an alternative proof, because a single function is not
sufficient to separate points in the two--dimensional space of the
parameters. 
\endexample

\subhead 6. Associated Sturm--Liouville operators \endsubhead
It is well known \footnote{One of the authors (G.~O.) is grateful to
F.~Alberto~Gr\"unbaum for drawing his attention to this fact.} that the
integral operator with the sine kernel, restricted to an arbitrary
finite interval, commutes with a certain Sturm--Liouville
differential operator
$$
Df=(pf')'+qf. 
$$
Specifically, take the interval $[-\tau,\tau]$ with $\tau>0$; then
$$
p(x)=x^2-\tau^2,\qquad q(x)=\pi^2x^2.
$$
Gaudin's proof of this fact, sketched in \cite{Me}, section 5.3, uses
a trick but the claim can also be verified by brute force.

Similar results also hold for the Airy kernel and the Bessel kernel,
and they turn out to be useful in the study of the corresponding
Fredholm determinants, see \cite{TW2, TW3}.

Now we shall produce analogous differential operators for the
kernels discussed in section 4 and for the Whittaker kernel.

\proclaim{Proposition 6.1} The integral operator with the $\sin/\sh$
kernel or the $\sh/\sh$ kernel, restricted to the interval
$[-\tau,\tau]$, $\tau>0$, commutes 
with the Sturm--Liouville operator $Df=(pf')'+qf$, where
$$
p(x)=\frac{\sh^2(Bx)-\sh^2(B\tau)}{B^2}\,,\qquad
q(x)=\frac{(B^2\pm A^2)\sh^2(Bx)}{B^2}
$$
and the plus sign is taken for the $\sin/\sh$ kernel while the minus
sign is taken for the $\sh/\sh$ kernel; the same is true in the limit
cases $A=0$ or $B=0$.
\endproclaim  

\proclaim{Proposition 6.2} The integral operator with the Whittaker
kernel restricted to the semi--infinite interval $(\tau,+\infty)$ with
arbitrary $\tau>0$, commutes with the Sturm--Liouville operator
$Df=(pf')'+qf$, where 
$$
\gather
p(x)=x(x-\tau)\,,\qquad
q(x)=-\,\frac{((a-x/2)^2-t)(x-\tau)}{x}\,,\\
a=(z+z')/2,\quad t=zz'.
\endgather
$$
\endproclaim  

\demo{Proof of the propositions} One checks by a direct computation
that the kernel satisfies the relation $D_xK(x,y)=D_yK(x,y)$ which
implies the desired claim. \qed
\enddemo

\subhead 7. Comparison with Poisson--Dirichlet \endsubhead
Like the processes $\PPz$, the Poisson--Dirichlet processes are closely
related to harmonic analysis on the infinite symmetric group. The
both families of processes play a similar role but on different levels.
Namely, the Poisson--Dirichlet processes describe the decomposition
on ergodic components for certain measures. Those measures live on a
compactification of the infinite symmetric group and are employed in
the construction of the generalized regular representations, see
\cite{KOV}. The processes $\PPz$ correspond to the dual level, as
they govern the decomposition of that representations. 

Our study makes it possible to compare the both families of
processes. 

The correlation functions of $\PDt$ are given by simple formulas,
those of $\PPz$ look much more complicated. This is already seen for
the first correlation functions.

Both $\PDt$ and $\PPz^+$ are simplified after lifting. But the
lifting of the former is a (non stationary) Poisson process while the
lifting of the latter is a less elementary object --- the fermion
process with the Whittaker kernel.

A similar conclusion can be made about the corresponding tail
processes. These are the standard Poisson process and a stationary
fermion process, respectively. This means, for instance, that the
asymptotic probability distribution of the ratio $\al_j/\al_{j+1}$
(as $j\to\infty$) looks quite differently.

Thus, the processes $\PPz$ seem to be much more sophisticated objects
than the Poisson--Dirichlet processes. 

On the other hand, as is shown in section 5 above, there is a rough
characteristic with respect to which $\PDt$ and $\PPz^+$ behave
similarly: in the both cases, the rate of decay of the $\al_j$'s is
that of a geometric progression.

\bigskip
\bigskip

{\smc A.~Borodin}: Department of Mathematics, The University of
Pennsylvania, Philadelphia, PA 19104-6395, U.S.A.  E-mail address:
{\tt borodine\@math.upenn.edu} 
\newline\indent
{\smc G.~Olshanski}: Dobrushin Mathematics Laboratory, Institute for
Problems of Information Transmission, Bolshoy Karetny 19, 101447
Moscow GSP-4, RUSSIA.  E-mail address: {\tt olsh\@ippi.ras.ru,
olsh\@glasnet.ru} 

\enddocument
\end